\documentclass[11pt,letterpaper]{amsart}
\usepackage[dvipsnames]{xcolor}
\usepackage{graphicx}
\usepackage{amssymb,amscd,amsthm,amsxtra}
\usepackage{latexsym}
\usepackage{epsfig}
\usepackage{amsmath}
\usepackage{bm}
\usepackage{amssymb}
\usepackage{amsthm}
\usepackage{xcolor}
\usepackage[linktocpage=true,colorlinks=true,linkcolor=Blue,citecolor=BrickRed,urlcolor=RoyalBlue]{hyperref}
\usepackage[alphabetic]{amsrefs}
\usepackage{comment}
\usepackage[T1]{fontenc}
\usepackage[utf8]{inputenc}

\usepackage[colorinlistoftodos,prependcaption,textsize=tiny]{todonotes}
\usepackage{caption}

\usepackage{enumitem} 
\usepackage[mathscr]{euscript} 

\usepackage{mathrsfs} 

\usepackage{mathtools}
\usepackage{xcolor}

\numberwithin{equation}{section} 

\def\dashint{\,\ThisStyle{\ensurestackMath{
            \stackinset{c}{.2\LMpt}{c}{.5\LMpt}{\SavedStyle-}{\SavedStyle\phantom{\int}}}
        \setbox0=\hbox{$\SavedStyle\int\,$}\kern-\wd0}\int}

\def \B {\mathcal{B}}

\newtheorem{theorem}{Theorem} 
\newtheorem{lemma}{Lemma}
\newtheorem{proposition}{Proposition}  

\newtheorem{definition}{Definition}

\numberwithin{equation}{section}

\title[Vázquez's SMP for linear and non linear PDEs with Unbounded Coefficients]{On the Vázquez Strong Maximum Principle for Linear and Nonlinear Elliptic Equations with Unbounded Coefficients and Optimal Conditions}

\author[C. S. Aguiar]{Christopher Silva Aguiar}
\address{Pontifical Catholic University of Rio de Janeiro -- PUC-Rio\\22451-900, G\'avea, Rio de Janeiro-RJ, Brazil}
\email{aguiarcsilva@gmail.com}

\makeatletter
\renewcommand{\enddoc@text}{}
\makeatother

\begin{document}



\begin{abstract}
In this paper, we present adaptations of the classical strong maximum principle for a certain class of equations, with particular emphasis on linear operators in non-divergence form and a new strategy to prove it. The novelty of our approach lies in the setting of the problem: the coefficients of the operators are unbounded, assumed to be merely integrable in the $L^{p,q}$-sense, and the source term $f$ on the right-hand side of the equation is treated under optimal integrability assumptions.
\end{abstract}


\maketitle


\section{Introduction}\label{Introduction}

The main objective of this article is to give a new proof to the Vázquez Strong Maximum Principle (VSMP) for equations of the type
 \begin{equation} \label{1}
  \mathcal{L}[u] := \textrm{tr}(A(x)M) + b(x)p + c(x)u  \leq  f(u)  \hspace{.5cm} \mbox{ in }  \Omega, 
\end{equation}

where $A$ is continuous and uniformly elliptic, and the lower-term coeficients are integrables. \\

The (SMPV) establishes that the Strong Maximum Principle (SMP), widely known in the literature, remains valid for nonnegative super-solutions of an equation of the type $(\ref{1})$, where the right-hand side term is non-homogeneous and nonlinear, governed by a function $f$. The study of these results began with Stampacchia in \cite{Prim} (but in a linear case, with $f(u) = c(x)u$) and the non linear case was later deepened by J. L. Vázquez in \cite{Vazquez}, where the following classical result was demonstrated:

\begin{theorem}[Vasques - 1984]\label{T1}
   Let $f : \mathbb{R} \longrightarrow \mathbb{R}$ be a continuous function such that: 
\begin{enumerate}
    \item $f \geq 0$, $f(0) = 0$;
    \item $f$ is nondecreasing on the interval $(0, \infty)$;
    \item $f$ satisfies the following condition:
    \begin{equation*}
        \int_{0}(F(s))^{-1/2}ds = \infty, \hspace{1cm} F(s) = \int_{0}^{s}f(t)dt.
    \end{equation*}
\end{enumerate}
Then, every classical super-solution to the equation
\begin{equation*}
    \Delta u = f(u) \hspace{0.2cm} \textrm{ in } \hspace{0.2cm} \Omega,
\end{equation*}
 satisfies the Strong Maximum Principle (SMP), that is $u=0$ or $u > 0$ in $\Omega$. 
\end{theorem}

 The proof of this theorem consists in build a sub-solution to compare it with the super-solution using the Comparison Principle (CP). This sub-solution was obtained from solutions of a Ordinary Differential Equation (ODE), that appears thanks to the change of variables and the boundness of the coefficients. This sub-solution is defined in a small annulus inside the initial domain and have a non-zero normal derivative on the boundary, getting a Hopf's Lemma for a radial problem. \\

The main focus of our work is around the condition $(3)$, where is known in the literature as the Vázquez Condition. While P. Bénilan, H. Brézis, and M. Crandall demonstrated the necessity of this hypothesis (see \cite{Sharp}), Vázquez shows that this rule about $f$ is sharp! 
Throughout this paper, we will refer the condition $(3)$ to it as $(VC)$.\\

Subsequently, several generalizations were obtained for more general operators, although still with limitation in the lower coefficients (bounded) and having the same strategy in our proofs. Some examples include:
\begin{itemize}
    \item In \cite{Lap1}, the authors proved the result for the Laplace operator augmented by a term dependent on the gradient's modulus of the solution, defined on all of $\mathbb{R}^N$, allowing the operator to have either a fully nonlinear or a quasilinear nature;
    \item In \cite{pqLap} was consider the $(p,q)$-laplacian operator;
    \item In the works \cite{Div1}, \cite{Div2}, and especially \cite{Div3}, the authors proved the result for general operators in divergence form, uniformly elliptic;
    \item For operators in non-divergence form, we cite \cite{ND}, where the authors established the result for the extremal Pucci operators with constant coefficients (covering, as a particular case, linear operators with bounded coefficients) for $f$ given in terms of the natural logarithm. However, the argument for $f$ according to the hypotheses described in Theorem 1 is analogous. We will state and prove this result later, as it will be necessary for our argument.
\end{itemize}

Up to the present moment, this result has only been obtained in the context of unbounded coefficients, with Sirakov and Souplet in your recent work \cite{BP}. They uses a different strategy to prove the result to a more general class of fully nonlinear operators, where encompassing the previous ones. The proof strategy relies on a contradiction argument. First, the existence of a subset $B \subset \Omega$ is shown, where the trace of the solution does not vanish on the boundary of $B$. Then, a another (Hölder-)continuous solution is constructed that coincides with the first one on $\partial B$. To obtain the contradiction, it is enough to show that this new solution is identically zero. To do so, a particular version of the (VSMP) is proved for smooth domains and continuous solutions, using a Harnack inequality with an optimized constant in its dependence on the coefficients. This result is then applied to the solution found, yielding the expected contradiction.

For other hand, they proved the theorem using the follow condition on $f$:  
    \begin{equation*}
        \limsup_{s \rightarrow 0} \frac{f(s)}{s(\ln s)^{2}} < \infty.
    \end{equation*}
Observe that this growth condition is more restrictive than the $(VC)$ present in Theorem \ref{T1} (where it can be seen taking $f(x) = x (ln x)^a (lnx)^b$, for $a,b \leq 2$). Isn't clear and obvious that the $(VC)$ was suitable and appropriate for the unlimited scenario, considering that it was purposed originally in the opposite context (bounded), making this issue remain open, even for the case envolving the linear operators.\\ 

Thus, the goal of this work is to solve one of this open problems, giving a new proof of the (VSMP), preserving the original and sharp condition $(VC)$, for the complete linear operator in non-divergence form. To do so, we adopt the following approach: we transform the original equation, whose low-order coefficients are unbounded, into a new equation in which these coefficients have better behavior (bounded). In this way, we can reuse the results already established in the literature.\\

The proof will be gived in two steps, were two essential tools is required. The first one is contained in the first step, is about the Zvonkin-type transformations. Widely used in Stochastic Partial Differential Equations (SPDEs), it was initially introduced in \cite{ZT}, whose the main goal is eliminates (or nullifyies) the first-order coefficient of the equation, called the "drift", thus allowing the recovery some fundamental results, such as the existence of classical solutions, and some other good properties. The second tool is the first eigenvalue along with the first eigenfunction (see, for example, \cite{AutV1}, \cite{AutV2}), which is relevant due to its applicability in recovering and generalizing important estimates seen in Partial Differential Equations (PDEs) theory, such as the (CP) and the Alexandroff-Bakeman-Pucci (ABP) estimates. \\

Now we will enunciate precisely our main result. To this, some important configurations is required.\\

Let $\Omega \subset \mathbb{R}^{n}$, $n \geq 2$ be a smooth and bounded domain. Consider a second-order linear operator in non-divergence form gived by 
\begin{equation}
    \mathcal{L}_{ND}[u] = \textrm{tr}(A(x)D^2 u ) + b(x) \cdot Du + c(x)u.
\end{equation}
We will always assume that $A(x) \in C(\Omega)$ satisfies: 
\begin{equation}
   \textrm{there exits constants }  0<\lambda \leq \Lambda \textrm{ such that }  \lambda I \leq A(x) \leq \Lambda I, \hspace{5mm} \forall x \in \Omega.
\end{equation}
The low-order coefficients are merely Lebesgue integrable, the minimum required for an (ABP) estimate to exist. Thus,
\begin{equation}
    b \in L^{q}(\Omega), c \in L^p(\Omega), \hspace{1mm} \textrm{ for } q > n \textrm{ and } p > p_{E}, 
\end{equation}
 where $p_{E} \in (n/2, n)$ represents the Euscariaza regime as seen in \cite{Exp}. 

Moreover, throughout this work, we will refer to $u \in C(\Omega)$ as a (sub-, super-)solution if $u$ satisfies the (in)equalities in the $L^q$-viscosity sense as in \cite{Caf}. \\

We now state our main result:

\begin{theorem}\label{T3}
Let $\Omega \subset  \mathbb{R}^n$, $n \geq 2$ and $b,f \in L^q(\Omega)$, $c \in L^p(\Omega)$, with $q>n$, $p > p_E \in (n/2,n)$ and $A \in C(\Omega)$ satisfying $\lambda I \leq A(x)\leq \Lambda I$, for some $0<\lambda \leq \Lambda$. Let $u$ being a $L^p-$viscosity solution to 
   \begin{equation*}
    \left \{ \begin{array}{rclrl} 
    \mathcal{L}_{ND}[u] \hspace{-.2cm} & \leq & \hspace{-.2cm} f(u) & \mbox{ in }  & \hspace{-.2cm} \Omega \\
     u(x) \hspace{-.2cm} & \geq & \hspace{-.2cm} 0 & \mbox{ in } & \hspace{-.2cm} \Omega\\
    u(x_{0}) \hspace{-.2cm} & = & \hspace{-.2cm} 0 & \mbox{ for some } & \hspace{-.2cm} x_0 \in \Omega. \end{array} \right.
\end{equation*}
with $f : \mathbb{R} \longrightarrow \mathbb{R}$ contínuous such that:

\begin{itemize}
    \item $f \geq 0$, $f(0) = 0$;
    \item $f$ is nondecreasing on the interval $(0, \infty)$;
    \item $f$ satisfies the $(VC)$ :
    \begin{equation*}
        \int_{0}(F(s))^{-1/2}ds = \infty, \hspace{1cm} F(s) = \int_{0}^{s}f(t)dt.
    \end{equation*}
\end{itemize}
Then $u \equiv 0$.
\end{theorem}

How said above, the proof will be gived in two steps. In the first one, the key is the Zvonkin Transform. We demonstrate the existence of this change of variables with the help of some tools from elliptic equations theory (Existence Theorem, ABP, and $C^{1,\alpha}$ regularity estimate), and this will give to us a new linear equation, whose first order coefficinent is null! Next, in second step, we use the first eigenvalue together with the first eigenfunction to convert our remaining equation with unbounded zero order coefficient into a new equivalent one, but that allow us to apply the (SMP) in a particular case. This result can be obtained with a simple adaptation of the proof of Lemma $5.3$ presented in the appendix in \cite{ND}, by substituting $f$ to satisfy the stated hypotheses.

Our paper is organized as follow: in section $2$ we established the necessary assumptions, correlate the essential definitions and we enunciate results known from the literature that will help us in the proof of our main result. In section $3$ we proceed with the proof of our main result.

 
\section{Preliminaries}\label{sec_prelim}

This section gathers the auxiliary tools that we use in our argumentation. Let's begin by defining and stating some results that we will use throughout the work, starting with the concept of the first eigenvalue and the first eigenfunction for an operator. All this structure  was provided in \cite{AutV1} and \cite{AutV2}. Consider $p<\infty$. 

\begin{definition} We say that the pair $\varphi, \kappa$ is an eigenfunction and an eigenvalue for a operator $F$ respectively, if $\varphi$ is a solution for
    \begin{equation*}
    \left \{ \begin{array}{rclrl} 
    F[\varphi]  \hspace{-.2cm} & = & \hspace{-.2cm} -\kappa \varphi_1 & \mbox{ in }  & \hspace{-.2cm} \Omega \\
     \varphi_1 \hspace{-.2cm} & > & \hspace{-.2cm} 0 & \mbox{ in } & \hspace{-.2cm} \Omega\\
    \varphi_1 \hspace{-.2cm} & = & \hspace{-.2cm} 0 & \mbox{ on } & \hspace{-.2cm} \partial \Omega. \end{array} \right.
\end{equation*}
\end{definition}

After setting the notions of eigenvalue and eigenfunction, we can see that this objects exists, in fact, for the operator $\mathcal{L}_{ND}$.

However, this notions will be important in the step two of our proof, in which, can be obtained in a more general setting. Consider $\mathcal{M}_{\lambda, \Lambda}^{\pm}$ represent the extremal Pucci operators given by:
\begin{equation*}
    \mathcal{M}_{\lambda, \Lambda}^{-}(M) = \inf_{\lambda I \leq A\leq \Lambda I } \textrm{tr}(AM) \hspace{2mm} \textrm{ and } \hspace{2mm}  \mathcal{M}_{\lambda, \Lambda}^{+}(M) = \sup_{\lambda I \leq A\leq \Lambda I } \textrm{tr}(AM).
\end{equation*}
These operators have some usefull properties, that we can see below:

\begin{lemma}\label{LEM}
Let $M, N \in S_n(\mathbb{R})$, $p, q \in \mathbb{R}^n$, and $\phi(x) \in C(\Omega)$ be such that $0 < a \leq \phi(x) \leq A$.
Then:
\begin{equation*}
    \mathcal{M}_{\lambda, \Lambda}^{-}(M) + \mathcal{M}_{\lambda, \Lambda}^{-}(N) \leq \mathcal{M}_{\lambda, \Lambda}^{-}(M+N) \leq \mathcal{M}_{\lambda, \Lambda}^{-}(M) + \mathcal{M}_{\lambda, \Lambda}^{+}(N),
\end{equation*}
\begin{equation*}
    \mathcal{M}_{\lambda a, \Lambda A}^{-}(M) \leq \mathcal{M}_{\lambda, \Lambda}^{-}(\phi M) \leq \mathcal{M}_{\lambda A, \Lambda a}^{-}(M),
\end{equation*}
\begin{equation*}
    -C(\lambda, \Lambda, n ) \vert p \vert \vert q \vert \leq \mathcal{M}_{\lambda, \Lambda}^{-}(p\otimes q) \leq \Lambda \vert p \vert \vert q \vert.
\end{equation*}
\end{lemma}

Finally, we can conclude the existence of eigenvalues and eigenfunctions for the classes of fully nonlinear operators depending on $\mathcal{M}_{\lambda, \Lambda}^{\pm}$, that was showed in \cite{KR}:

\begin{theorem}\label{EAU2} 
Let $n<q < \infty$ and $c \in L^p(\Omega)$, with $c(x) \leq 0$ for all $x$. For $E_q := W_{loc}^{2,q}(\Omega) \cap C(\overline{\Omega})$, define
\begin{equation*}
    \Psi_0^+(\mathcal{L}^-,\Omega,\lambda) = \{ \psi \in E_q \hspace{1mm}; \hspace{1mm} \mathcal{L}^-[\psi] + \lambda \psi \leq 0  \hspace{1mm} in \hspace{1mm}\Omega \hspace{1mm} and\hspace{1mm} \psi > 0 \hspace{1mm} in \hspace{1mm} \Omega  \}
\end{equation*}
and
\begin{equation*}
    \Psi_0^-(\mathcal{L}^-,\Omega,\lambda) = \{ \psi \in E_q \hspace{1mm}; \hspace{1mm} \mathcal{L}^-[\psi] + \lambda \psi \geq 0  \hspace{1mm} in \hspace{1mm}\Omega \hspace{1mm} and \hspace{1mm} \psi < 0 \hspace{1mm} in \hspace{1mm} \Omega \},
\end{equation*}
where
\begin{equation*}
    \mathcal{L}^-[u] = \mathcal{M}_{\lambda,\Lambda}^{-}(D^2u) + c(x)u.
\end{equation*}
There are functions $\varphi_1^+, \varphi_1^- \in E_p$ such that if we define
\begin{equation*}
    \lambda_1^+(F, \Omega) := \sup \{\lambda \hspace{1mm} ; \hspace{1mm} \Psi_0^+(F,\Omega,\lambda) \neq \emptyset \}
\end{equation*}
and
\begin{equation*}
    \lambda_1^-(F, \Omega) := \sup \{\lambda \hspace{1mm} ; \hspace{1mm} \Psi_0^-(F,\Omega,\lambda) \neq \emptyset \}
\end{equation*}
then each pair $(\lambda_1^{+}, \varphi_1^{+})$ and $(\lambda_1^{-}, \varphi_1^{-})$ satisfies
  \begin{equation*}
\begin{minipage}{0.40\textwidth}
  \begin{align*}
      \left \{ \begin{array}{rclrl} 
    \mathcal{L}^-[\varphi_1^-]  \hspace{-.2cm} & = & \hspace{-.2cm} -\lambda_1^- \varphi_1^- & \mbox{ in }  & \hspace{-.2cm} \Omega \\
     \varphi_1^- \hspace{-.2cm} & < & \hspace{-.2cm} 0 & \mbox{ in } & \hspace{-.2cm} \Omega\\
    \varphi_1^- \hspace{-.2cm} & = & \hspace{-.2cm} 0 & \mbox{ on } & \hspace{-.2cm} \partial \Omega \end{array} \right.
    \end{align*}
\end{minipage} \textrm{ and }
\begin{minipage}{0.40\textwidth}
    \begin{align*}
      \left \{ \begin{array}{rclrl} 
    \mathcal{L}^-[\varphi_1^+]  \hspace{-.2cm} & = & \hspace{-.2cm} -\lambda_1^+ \varphi_1^+ & \mbox{ in }  & \hspace{-.2cm} \Omega \\
     \varphi_1^+ \hspace{-.2cm} & > & \hspace{-.2cm} 0 & \mbox{ in } & \hspace{-.2cm} \Omega\\
    \varphi_1^+ \hspace{-.2cm} & = & \hspace{-.2cm} 0 & \mbox{ on } & \hspace{-.2cm} \partial \Omega. \end{array} \right.
    \end{align*}
\end{minipage}
\end{equation*}    
\end{theorem}

Notice that in both cases, we can take $\varphi$ normalized in $x_0$, that is, $\varphi_1(x_0) = \varphi_1^+(x_0) = \varphi_1^-(x_0) = 1$ and there exists a constant $C$ depending only on the initial data of the problems ($x_0, \Omega, \lambda, \Lambda, b, c$) such that $\varphi_1(x_0),\varphi_1^+(x_0), \varphi_1^- < C$. 

Moreover, we are showing our result in the extremal case. This means that, for any uniformly elliptical operator $F$ the result still valid as long as $F$ is convex, for example, the other extremal Pucci operator $\mathcal{L}^+$ and $\mathcal{L}_{ND}$. \\

Now, we will state some auxiliary results from the classical theory of (PDEs) in which we will use throughout this paper.

\begin{theorem}[ABP Estimate]\label{L1}
If diam$(\Omega) \leq 1$, $\|b\|_{L^{q}(\Omega)}\leq 1$, $f,c \in L^p(\Omega)$, with $q>n$, $p >p_E$ and $w \in C(\Omega)$ is a supersolution for
\begin{equation*}
    \left \{ \begin{array}{rclrl} 
    \mathcal{M}_{\lambda, \Lambda}^{-}(D^2w) - b(x)\vert Dw \vert  \hspace{-.2cm} & \leq & \hspace{-.2cm} f(x) & \mbox{ in }  & \hspace{-.2cm} \Omega \\
     w(x) \hspace{-.2cm} & \geq & \hspace{-.2cm} 0 & \mbox{ on } & \hspace{-.2cm} \partial\Omega.\\
     \end{array} \right.     
\end{equation*}
Then for some constant $C = C(n,q, \lambda, \Lambda )$, we have that
\begin{equation*}
    \sup_{\Omega}w \leq C\|f\|_{L^{n}(\Omega^+)},
\end{equation*}
    where $\Omega^+ = \{x \in \Omega \hspace{1mm};\hspace{1mm} w(x) > 0\}$.
\end{theorem}

\begin{proof}
    Apply the Theorem $2.9$ of \cite{ABP} in each connected component of $\Omega^+.$
\end{proof}

\begin{theorem}[Existence and Regularity]\label{ER}
Let $\Omega$ satisfy a uniform exterior cone condition. Then there exists a unique $L^p$-viscosity solution $u \in W_{loc}^{2,p}(\Omega) \cap C(\overline{\Omega})$ of
\begin{equation*}
     \left \{ \begin{array}{rclrl} 
    \textrm{tr}(A(x)D^2(u)) + b(x)Du \hspace{-.2cm} & = & \hspace{-.2cm} f(x) & \mbox{ in }  & \hspace{-.2cm} \Omega \\
     u(x) \hspace{-.2cm} & = & \hspace{-.2cm} g(x) & \mbox{ on } & \hspace{-.2cm} \partial\Omega,\\
     \end{array} \right.    
\end{equation*}
where $b,f \in L^p(\Omega)$ and $g \in C(\partial \Omega)$. Moreover, for every $\Omega' \subset \overline{\Omega'} \subset \Omega$, there exists a constant $C = C(n,p, \lambda, \Lambda, \omega_A, \|b\|_{L^p(\Omega)}, diam(\Omega), dist (\Omega', \partial \Omega))$ such that
\begin{equation*}
    \|u\|_{W^{2,p}(\Omega')} \leq C \left(\|u\|_{L^{\infty}(\Omega)} + \|f\|_{L^p(\Omega)} + \|g\|_{C(\partial \Omega)}\right).
\end{equation*}
\end{theorem}

\begin{proof}
    Just apply the Corollary $3.10$ and Lemma 2.5 in \cite{Caf} for the operator above.
\end{proof}
It's important to remember that using the Sobolev's Embbedings we have $C^{1, \alpha}$ estimates for some $\alpha = \alpha(p)$.

\section{Bounded Case}
To demonstrate the Strong Maximum Principle of Vázquez in the unbounded setting, we must first establish it in the bounded case. Therefore, the goal of this section is to present that proof. We begin with the main tool of the argument: the construction of a special barrier function.

\begin{proposition}
  Let $k_1,k_2,r_1,b$ and $v_1$ be positive constants. Then the boundary value problem
      \begin{equation*}
    \left \{ \begin{array}{rclrl} 
    v'' \hspace{-.2cm} & = & \hspace{-.2cm} \displaystyle M\left(k_1m(v') + b \vert v' \vert+ k_2g(v) \right), & \mbox{ for }  & \hspace{-.2cm} 0 < r < r_1,  \\
     v(0) \hspace{-.2cm} & = & \hspace{-.2cm} 0, & & \hspace{-.2cm} \\
    v(r_1) \hspace{-.2cm} & = & \hspace{-.2cm} v_1 .& & \hspace{-.2cm} \end{array} \right.
\end{equation*}
admits a $C^1$ solution defined on $[0,r_1]$ such that $v,v' \geq 0$, and if $g$ satisfies the $(VC)$ condition, then $0 < v < v_1$ in $(0,r_1)$.
\end{proposition}

\begin{proof}
    Let us divide the topics to clarify where and how each item is demonstrated.\\
    
    1- Existence.\\

For this case, it suffices to observe that an ODE is essentially a PDE in one dimension. Therefore, the method of sub- and supersolutions (Perron's method) applies, together with the comparison principle. In this way, by taking $w_1 = v_1$  as a supersolution and $w_2 = 0$ as a subsolution, we obtain the existence of a solution satisfying $w_1 \geq v \geq 0$. Moreover, $v$ is of class $C^1$ (in fact, even smoother), since the right-hand side of our equation is continuous. \\

2- $v' \geq 0$.\\

Suppose, by contradiction, that this does not occur. Then there exists a point $r_0$ such that $v'(r_0) < 0$, and consequently, there exists an open interval $I$ where $v'(r) < 0$ for all $r \in I$. Define the auxiliary function $\varphi(r) = e^{-k_1\frac{\Lambda}{\lambda}r}v'(r)$. Then, by differentiating in some $r \in I$, using the differential equation, and the definition of $M$ we obtain:
\begin{equation*}
     \begin{array}{rclrl} 
\varphi'(r) \hspace{-.2cm} & = & \displaystyle \hspace{-.2cm} - \frac{k_1\Lambda}{\lambda}e^{\frac{k_1\Lambda}{\lambda}r}v'(r) + e^{\frac{k_1\Lambda}{\lambda}r}v''(r) \vspace{.2cm} \\
      & = &  \hspace{-.2cm} \displaystyle - \frac{k_1\Lambda}{\lambda}e^{\frac{k_1\Lambda}{\lambda}r}v'(r) + e^{\frac{k_1\Lambda}{\lambda}r}\displaystyle M\left(k_1m(v') + b \vert v' \vert+ k_2g(v) \right) \vspace{.2cm} \\
      & = & \hspace{-.2cm} \displaystyle - \frac{k_1\Lambda}{\lambda}e^{-k_1r}v'(r)+ e^{-k_1r}\frac{ k_1\Lambda v'}{\lambda} + e^{-k_1r}\frac{\left(b \vert v' \vert+ k_2g(v)\right)}{\lambda} \vspace{.2cm}\\
      & = & \hspace{-.2cm} \displaystyle \frac{e^{-k_1r}}{\lambda}\left(b \vert v' \vert+ k_2g(v)\right),
    \end{array} 
\end{equation*}
concluding that $\varphi'(r) \geq 0$ on $I$, which implies that $\varphi$ is non-decreasing. Therefore, $v'(r) \geq 0$, which is a contradiction, since we assumed $v' < 0$ on $I$.

3- $0<v<v_1$.\\

Now suppose that Vázquez's condition holds. Let $r_0$ the maximum value for which $v(r) = 0$. Thus, by applying a change of variables, we obtain that
\begin{equation}\label{EQ55}
    \int_{r_0}^{r_1} \frac{v'(r)}{\sqrt{F(v(r))}} dr = \int_{0}^{v_1}\frac{1}{\sqrt{F(s)}}ds = \infty.
\end{equation}
By other hand, take $w = (v')^2$. Note that, as we proved in the previous item, we have that $v'\geq 0$, then the functions $M$ and $m$ takes a concrete value, making the equation satisfied by $v$ turn into
\begin{equation*}
     v''  =  k_1\frac{\lambda}{\Lambda}v' + \frac{b}{\Lambda} \vert v' \vert+ \frac{k_2}{\Lambda}g(v)
\end{equation*}
Then, applying the fundamental theorem of calculus and the chain rule, along with the equation, we obtain that
\begin{equation*}
    \frac{k_2}{\Lambda}F'(v) = \frac{k_2}{\Lambda}g(v)\cdot v' =  \left(v'' - \frac{(k_1\lambda + b)}{\Lambda}v'\right) \cdot v'
\end{equation*}
Therefore, defining another auxiliary function given by
\begin{equation*}
    h(r) = e^{-2\frac{(k_1\lambda+b)}{\Lambda}r}w(r),
\end{equation*}
we have
 \begin{equation*}
     \begin{array}{rclrl} 
    \left(e^{-2\frac{(k_1\lambda+b)}{\Lambda}r}w(r)\right)' \hspace{-.2cm} & = & \displaystyle \hspace{-.2cm} -2\frac{(k_1\lambda+b)}{\Lambda}e^{-2\frac{(k_1\lambda+b)}{\Lambda}r}w(r) + e^{-2\frac{(k_1\lambda+b)}{\Lambda}r} w'(r) \vspace{.2cm} \\
      & = & \displaystyle \hspace{-.2cm} -2\frac{(k_1\lambda+b)}{\Lambda}e^{-2\frac{(k_1\lambda+b)}{\Lambda}r}v'(r)\cdot v'(r) + e^{-2\frac{(k_1\lambda+b)}{\Lambda}r} 2v''(r) \cdot v'(r) \vspace{.2cm} \\
      & = & \displaystyle \hspace{-.2cm} 2e^{-2\frac{(k_1\lambda+b)}{\Lambda}r}v'(r)\cdot\left( -\frac{(k_1\lambda+b)}{\Lambda} v'(r) + v''(r)\right) \vspace{.2cm} \\
      & = & 2\displaystyle \frac{k_2}{\Lambda}F'(v)e^{-2\frac{(k_1\lambda+b)}{\Lambda}r}.
    \end{array} 
\end{equation*}
Now, assume that $v'(r_0) = 0$. In this way, integrating both sides from $r_0$ to $r$ and using the fact that the function $\frac{1}{e^x}$ is decreasing, this gives us that: on the one hand,
\begin{equation*}
    \int_{r_0}^{r}\left(e^{-2\frac{(k_1\lambda+b)}{\Lambda}r}w(r)\right)'dr = e^{-2\frac{(k_1\lambda+b)}{\Lambda}r}w(r) \geq e^{-2\frac{(k_1\lambda+b)}{\Lambda}r_1}w(r),
\end{equation*}
and, on the other hand, 
\begin{equation*}
   \int_{r_0}^{r} 2\displaystyle \frac{k_2}{\Lambda}F'(v)e^{-2\frac{(k_1\lambda+b)}{\Lambda}r}dr \leq  2\displaystyle \frac{k_2}{\Lambda}e^{-2\frac{(k_1\lambda+b)}{\Lambda}r_0}\int_{r_0}^{r}F'(v)dr = 2\displaystyle \frac{k_2}{\Lambda}e^{-2\frac{(k_1\lambda+b)}{\Lambda}r_0} F(v(r)).
\end{equation*}
Therefore,
\begin{equation*}
    2\displaystyle \frac{k_2}{\Lambda}e^{-2\frac{(k_1\lambda+b)}{\Lambda}r_0}F(v(r)) \geq e^{-2\frac{(k_1\lambda+b)}{\Lambda}r_1}w(r),
\end{equation*}
Where, by manipulating the terms, we obtain that
\begin{equation*}
    \frac{w(r)}{F(r)} \leq  2\displaystyle \frac{k_2}{\Lambda} \frac{e^{-2\frac{(k_1\lambda+b)}{\Lambda}r_0}}{e^{-2\frac{(k_1\lambda+b)}{\Lambda}r_1}} = 2\displaystyle \frac{k_2}{\Lambda} e^{2\frac{(k_1\lambda+b)}{\Lambda}[r_1 - r_0]} \Longleftrightarrow   
    \frac{v'(r)}{\sqrt{F(r)}} \leq  \sqrt{2\displaystyle \frac{k_2}{\Lambda}}e^{\frac{(k_1\lambda+b)}{\Lambda}[r_1 - r_0]}
\end{equation*}
Consequently, integrating from $r_0$ to $r_1$ we get
\begin{equation*}
    \int_{r_0}^{r_1}\frac{v'(r)}{\sqrt{F(v(r))}} dr \leq \int_{r_0}^{r_1}\sqrt{2\displaystyle \frac{k_2}{\Lambda}}e^{\frac{(k_1\lambda+b)}{\Lambda}[r_1 - r_0]} = \sqrt{2\displaystyle \frac{k_2}{\Lambda}}e^{\frac{(k_1\lambda+b)}{\Lambda}[r_1 - r_0]} [r_1 - r_0] < \infty.
\end{equation*}
This is a contradiction with $(\ref{EQ55})$. Thus $v'(r_0) > 0$, wich implies that $r_0 = 0$.
\end{proof}
We can give now the proof of the Strong Maximum principle in a bounded setting.

\begin{theorem}\label{VBounded}
    Let $u \in C(\B_1)$ be a viscosity super-solution to 
  \begin{equation*}
    \left \{ \begin{array}{rclrl} 
    \mathcal{M}_{\lambda, \Lambda}^{-}(D^2u) - b\vert Du \vert - c u  \hspace{-.2cm} & \leq & \hspace{-.2cm} f(u) & \mbox{ in }  & \hspace{-.2cm} \B_1 \\
     u(x) \hspace{-.2cm} & \geq & \hspace{-.2cm} 0 & \mbox{ in } & \hspace{-.2cm} \B_1 \end{array} \right.
\end{equation*}

Then, either $u \equiv 0$ or $u>0$ in $\Omega$.
\end{theorem}
\begin{proof}
Let us consider a particular case, taking $u \in C^1$, which will be of great use to us.
The proof proceeds by contradiction. Suppose the result is false. Then the sets
\begin{equation*}
    \Omega^+ = \{x \in \B_1 \,;\, u(x) > 0\} \,\,\ \textrm{ e } \Omega^0 = \{x \in \B_1 \,;\, u(x) = 0\}
\end{equation*}
are both non-empty. Furthermore, since $u$ is continuous, $\Omega^+$ is a open set. Therefore, it is possible to choose a point $\Tilde{x} \in \Omega^+$ such that $d(\Tilde{x}, \partial \Omega^0) < d(\Tilde{x}, \partial \B_1)$ and there exists a ball $\B_R(\Tilde{x})$ such that $\partial\B_R(\Tilde{x}) \cap \partial\Omega^0 \neq 0$.  Let $x_0$ be a point in this intersection. Thus, $u(x_0) = 0$, and we may assume, possibly after reducing the radius $R$, that $u<\delta$ in $\B_R(\Tilde{x})$. Finally, we define the annulus $E_R = \B_{R}(\Tilde{x}) \setminus \overline{\B}_{\frac{R}{2}}(\Tilde{x})$ and the quantity $\mu = \min_{\partial\B_{R/2}(\Tilde{x})}u < \delta$. 

The main idea is to construct a barrier function to extract information about the derivative of $u$ at $x_0$, ultimately leading to a contradiction. This function will be constructed through a duality between this PDE and an associated ODE. First, note that we may assume $c = 0$; otherwise, define $g(t) = f(t) + ct$, so that if $u$ s a solution of the original equation, it also solves:
 \begin{equation*}
    \left \{ \begin{array}{rclrl} 
    \mathcal{M}_{\lambda, \Lambda}^{-}(D^2u) - b\vert Du \vert   \hspace{-.2cm} & \leq & \hspace{-.2cm} g(u) & \mbox{ in }  & \hspace{-.2cm} \B_1 \\
     u(x) \hspace{-.2cm} & \geq & \hspace{-.2cm} 0 & \mbox{ in } & \hspace{-.2cm} \B_1 \end{array} \right.
\end{equation*}
Moreover, $g$ atisfies all the same assumptions as $f$, including the $(VC)$, Since the square root is a subadditive function. According to the previous proposition, let $\hat{u}(x) = v(R-\vert x-\Tilde{x}\vert, k_1, 1, R/2,\mu)$, where $v$ solves the following problem::
      \begin{equation*}
    \left \{ \begin{array}{rclrl} 
    v'' \hspace{-.2cm} & = & \hspace{-.2cm} \displaystyle M\left(-\frac{N-1}{r}m(v') + b \vert v' \vert+ g(v) \right), & \mbox{ for }  & \hspace{-.2cm} 0 < r < R/2,  \\
     v(0) \hspace{-.2cm} & = & \hspace{-.2cm} 0, & & \hspace{-.2cm} \\
    v(R/2) \hspace{-.2cm} & = & \hspace{-.2cm} \mu. & & \hspace{-.2cm} \end{array} \right.
\end{equation*}
It is known that choosing a adequade $k_1$ the function $\hat{u}$ satisfies:
 \begin{equation*}
  \begin{array}{rclrl} 
    \mathcal{M}_{\lambda, \Lambda}^{-}(D^2u) - b\vert Du \vert   \hspace{-.2cm} & \geq & \hspace{-.2cm} g(u) & \mbox{ in }  & \hspace{-.2cm} E_R. \\
    \end{array} 
\end{equation*}
We now prove that $u \geq \hat{u}$ in $E_R$. Suppose, for contradiction, that this is not the case. Then the open set $\mathcal{O} = E_R \cap \{ u < \hat{u} \}$ is non-empty. By construction and using the monotonicity of $g$, we have:
 \begin{equation*}
     \left \{\begin{array}{rclrl} 
   \mathcal{L}^-[\hat{u}]-\mathcal{L}^-[u]& \geq &  \hspace{-.2cm} g(\hat{u}) - g(u) \geq 0   & \mbox{ in }  & \hspace{-.2cm} E_R \\
   \hat{u}& \leq &u& \mbox{ on }  & \hspace{-.2cm} \partial\mathcal{O} \subset \partial E_R \cup \{\hat{u} = u\}
    \end{array} \right.
\end{equation*}
Therefore, by the comparison principle, $\hat{u} \leq u$ in $\mathcal{O}$, contradicting the definition of $\mathcal{O}$. Hence, $u \geq \hat{u}$ in $E_R$. We can now use this fact, combined with $\hat{u}'(0) > 0$, to observe that 
\begin{equation*}
    \liminf_{h\rightarrow 0} \frac{1}{h}[u(x_0 + h(\Tilde{x} - x_0)) - u(x_0)] \geq  \liminf_{h\rightarrow 0} \frac{1}{h}[\hat{u}(x_0 + h(\Tilde{x} - x_0)) - \hat{u}(x_0)] > 0.
\end{equation*}
However, since $u \in C^1$, this implies $\nabla u(x_0) > 0$, contradicting the fact that $\nabla u(x_0) = 0$ because $x_0$  is a local minimum of $u$. Therefore, either $u>0$ or $u=0$ in $\B_1$. \\

Now, to the general case! Assume that $u$ is not identically zero. Then, by continuity, there exists a small ball $\B_R(\hat{x}) \subset \Omega$ such that the graph of the function $\hat{u}$ provided by the previous proposition, is strictly positive on the boundary $\partial B_R(\tilde{x})$ and hence $\beta = \min_{\partial \B_R(\tilde{x})}\hat{u}$ is well-defined. Now define a new function $v$ as the solution to the problem:
 \begin{equation*}
     \left \{\begin{array}{rclrl} 
   -\mathcal{L}^-[v] \hspace{-.2cm}& = &  \hspace{-.2cm} - g(v)    & \mbox{ in }  & \hspace{-.2cm} B_R(\tilde{x}) \\
   v& = & \min\{\delta/2,\beta \}  & \mbox{ on }  & \partial B_R(\tilde{x})
    \end{array} \right.
\end{equation*}
It is known from standard theory that the above problem admits a unique $C^{1,\alpha}$ solution. Moreover, by the construction method, we obtain $0\leq v \leq \delta/2$. Applying the previous $C^1$ case to $v$, we deduce that $v > 0$ in $\B_R(\tilde{x})$. Then, If we now show that $u \geq v$ em $\B_R(\tilde{x})$, we finish. But this is an immediate consequence of the construction and the comparison principle, which proves the positivity of $u$ in some ball inside the domain. Having established the result for one interior ball, we can extend it to any ball intersecting $\B_R(\tilde{x})$ and then to the entire domain via a covering argument. Begin by covering the domain with balls around and intersecting the initial ball, and then extract a finite subcover to prove positivity on any compact subset of the domain.

\end{proof}

\section{Proof of Main Result}
First, let us reduce the problem. Observe that it is sufficient to prove the results for the case where $\Omega =\B_1$. Instead, fix $x_0 \in \Omega$. Since $\Omega$ is an open, bounded, and regular set, we can find an $r>0$ and a point $y_0 \in \Omega$ such that $x_0 \in \B_r(y_0) \subset \Omega$. Therefore, by performing the change of variables $z=\frac{x-y_0}{r}$, we have that $u(x)=v(z)$ solves a new equation with the same structure as the previous one, but with the domain now being $\B_1$ and with different ellipticity constants. From this point, it is enough to apply the theorem that will be demonstrated. Moreover, we can assume that $x_0 = 0$, using the translation $v = u-u(x_0)$.\\

The key to proving this result lies in the following argument: let $\Psi$ be a change of variables function to be determined, and let $y=\Psi(x)$. Upon applying the original equation, we obtain a function $v(y):=u(x)$, which solves:
\begin{equation*}
    \textrm{tr}(D\Psi^{T}A(x)D\Psi(x)D^2v(y)) + \overline{b}(y)Dv(y) + \overline{c}(y)v(y) \leq f(v(y)),
\end{equation*}
where
\begin{equation*}
    \overline{b}(y) = (\overline{b}_1(y), \ldots, \overline{b}_n(y)), \textrm{ com } \overline{b}_k(y) = \textrm{tr}(A(x) D^2\Psi_k(x))  +  b(x) D(\Psi_k(x)) 
\end{equation*}
and
\begin{equation*}
 \overline{c}(y) = c(x).
\end{equation*}
Our goal is to annul the coeficients $\overline{b}$ and $\overline{c}$. To achieve this, we will choose the change of variables as a small perturbation of the identity, small enough to keep the bijectivity.

\begin{proposition}\label{Prop1}
    In the same conditions of Theorem \ref{T3}, there exists a bijective change of variables $\Psi : \Omega \longrightarrow \Omega$ such that if $u$ solves 
    \begin{equation*}
    \left \{ \begin{array}{rclrl} 
    \mathcal{L}_{ND}[u] \hspace{-.2cm} & \leq & \hspace{-.2cm} f(u) & \mbox{ in }  & \hspace{-.2cm} \B_1 \\
     u(x) \hspace{-.2cm} & \geq & \hspace{-.2cm} 0 & \mbox{ in } & \hspace{-.2cm} \B_1\\
    u(0) \hspace{-.2cm} & = & \hspace{-.2cm} 0 &  &  \end{array} \right.
\end{equation*}
then $v(y):=u(x)$ solves
\begin{equation*}
    \left \{ \begin{array}{rclrl} 
    \textrm{tr}(D\Psi^{T}A(x)D\Psi(x)D^2v(y)) + \overline{c}(y)v(y) \hspace{-.2cm} & \leq & \hspace{-.2cm} f(v(y)) & \mbox{ in }  & \hspace{-.2cm} \Psi(\B_{r_0}) \\
     v(y) \hspace{-.2cm} & \geq & \hspace{-.2cm} 0 & \mbox{ in } & \hspace{-.2cm} \Psi (\B_{r_0})\\
    v(y_{0}) \hspace{-.2cm} & = & \hspace{-.2cm} 0 & \mbox{ for } & \hspace{-.2cm} y_0 = \Psi(0), \end{array} \right.
\end{equation*}
for some $r_0 > 0.$
\end{proposition}
\begin{proof}
     We are looking for a function $S(x)$ and a number $r_0 > 0$ for
\begin{equation*}
    \Psi(x) := x + S(x) = (x_1 + S_1(x), \ldots, x_n + S_n(x))
\end{equation*}
be bijective in a neighborhood $\B_{r_0}$ of the origin ($r_0$ to be found) and $\overline{b}$ be null. Then, applying the change of variables in definition of $\overline{b}_k$, we have:
\begin{equation*}
\begin{array}{rclrl}
 \overline{b}_k(y)  & \hspace{-.2cm} = & \hspace{-.2cm}\textrm{tr}(A(x) [D^2(x_k + S_k(x))]) + b(x)(D(x_k + S_k(x))) \\ 
\ & \hspace{-.2cm} = & \hspace{-.2cm}  \textrm{tr}(A(x)D^2S_k(x)) + b(x)DS_k(x)  + b_k(x).
\end{array}
\end{equation*}

This means that we can choose each $S_k(x)$ being a solution of the equation
 \begin{equation*}
    \Bigg \{ \begin{array}{rclrl} 
     \textrm{tr}(A(x)D^2S_k(x)) + b(x)DS_k(x) \hspace{-.2cm} & = & \hspace{-.2cm} -b_k(x) & \mbox{ in }  & \hspace{-.2cm} \B_r \\
     S_k(x) \hspace{-.2cm} & = & \hspace{-.2cm} 0 & \mbox{ in } & \hspace{-.2cm} \partial \B_r, \end{array} 
\end{equation*}
For some initial $0<r<1$. First, observe that such $S_k$ exists, is unique, and have $C^{1,\alpha}$ estimates by the Theorems \ref{ER} and \ref{L1}.

Now, taking another change of variables gived by $z = x/r$, we gets a new function $T_k(z) := S_k(x)$ that solves
\begin{equation*}
    \Bigg \{ \begin{array}{rclrl} 
    \textrm{tr}(\Tilde{A}(z)D^2T_k(rz)) + rb(rz)DT_k(x) \hspace{-.2cm} & = & \hspace{-.2cm} -r^2b(rz) & \mbox{ in }  & \hspace{-.2cm} \B_1 \\
     T_k(x) \hspace{-.2cm} & = & \hspace{-.2cm} 0 & \mbox{ in } & \hspace{-.2cm} \partial \B_1. \end{array} 
\end{equation*}

 Now, observe that the lower-order coefficient and the right-hand side depend on $r$. This will be useful because, when evaluating the $L^p$-norm of the term $rb(rz)$, we obtain:
 
\begin{equation}\label{Aux4}
    \| r b(rz) \|_{L^{q}(\B_1)} = r^{1 - \frac{n}{q}}  \| b \|_{L^{q}(\B_r)} \leq B r^{1 - \frac{n}{q}}  \underset{r \rightarrow 0}{\longrightarrow} 0,
\end{equation}

since $q>n$ and taking $B = \| b \|_{L^{q}(\B_r)} < \infty $. Applying the (ABP) estimates (Theorem \ref{L1}) we obtain

\begin{equation}\label{Aux5}
    \sup_{\B_1} \vert T_k \vert \leq C \|  r^2 b(rz) \|_{L^{p}(\B_1)} \leq  C r^{2 - \frac{n}{p}}  \underset{r \rightarrow 0}{\longrightarrow} 0,
\end{equation}
 where we used that $p > p_E$ and $p_E \in (n/2, n)$, which implies $\frac{n}{p} < \frac{n}{p_E} < 2$. Therefore, by the $C^{1, \alpha}$ estimates we have: 
\begin{equation*}
    \sup_{\B_r} \vert DS_k \vert \leq C(\|S_k\|_{L^{\infty}(B_1)} +\|  r^2 b(rz) \|_{L^{p}(\B_1)}),
\end{equation*}
and combining this with $(\ref{Aux4})$ and $(\ref{Aux5})$, we conclude that 
\begin{equation*}
    \sup_{B_r} \vert DS_k \vert \underset{r \rightarrow 0}{\longrightarrow} 0.
\end{equation*}
So, there exists $r_0 > 0$ such that, for $r<r_0$, 
\begin{equation*}
    \sup_{\B_r} \vert DT_k \vert \leq \frac{1}{2}.
\end{equation*}
Recall that, by the definition of $T$ and $S$, we have: 
\begin{equation*}
    \sup_{\B_1} \vert T_k \vert = \sup_{\B_r} \vert S_k \vert \hspace{0.2cm} \textrm{ e } \hspace{0.2cm} \sup_{\B_1} \vert DT_k \vert = \sup_{\B_r} \vert DS_k \vert.
\end{equation*}
Thus, if $r<r_0$, we have that $D\Psi = I + DS$ is bounded, that is, 
\begin{equation*}
    \frac{1}{2} I \leq D\Psi \leq \frac{3}{2} I
\end{equation*}
and therefore, is invertible for every $x \in \B_r$. By the Inverse Function Theorem, $\Psi$ is a bijection in this neighborhood. 

\end{proof}

Observe that, after this result, that is, after applying the change of variables $\Psi$, we remais with a solution $v$ of: 
\begin{equation*}
    \left \{ \begin{array}{rclrl} 
    \textrm{tr}(\overline{A}(y)D^2v(y)) + \overline{c}(y)v(y)  \hspace{-.2cm} & \leq & \hspace{-.2cm} f(v(y)) & \mbox{ in }  & \hspace{-.2cm} B \\
     v(y) \hspace{-.2cm} & \geq & \hspace{-.2cm} 0 & \mbox{ in } & \hspace{-.2cm} B\\
    v(y_0) \hspace{-.2cm} & = & \hspace{-.2cm} 0 & \mbox{ where } & \hspace{-.2cm} y_0 = \Psi(0). \end{array} \right.
\end{equation*}
where $B = \Psi(\B_{r_0})$ is a open and regular neighborhood of $y_0$.

For the next step, we will prove a more general case, considerating fully nonlinear operators. To this, let us consider the operator $G$ to be defined $S_n(\mathbb{R}) \times \mathbb{R} \times \Omega$. Moreover we will assume:\\
 \vspace{-3mm}

 \begin{equation*}
  \hspace{-4cm} (H_0) \hspace{3.8cm}   G(M,0,x) \textrm{ is contínuous in } S_n(\mathbb{R}) \times \Omega;
 \end{equation*}
\vspace{-1mm}
 \begin{equation*}
 \hspace{-4.5cm}  (H_1) \hspace{4cm} G(M,0,x) \textrm{ is convex in } S_n(\mathbb{R}) \times \mathbb{R}.
 \end{equation*}\\
\vspace{-0.5mm}
And finally, let us assume that $G$ to be $1-$homogeneos and uniformly elliptic: for some $\delta > 0$ and for all $M,N \in S_n(\mathbb{R})$, $u,v \in \mathbb{R}$, $x \in \Omega$ and $t\geq0$, we have

\begin{equation*}
  \hspace{-4.5cm}  (H_2) \hspace{5.2cm} G(tM,tu,x) = tG(M,u,x)
\end{equation*}
and\\
    \begin{equation*}
\begin{array}{rclrl} 
 \hspace{-2cm} (D_F) \hspace{3cm} \mathcal{M}_{\lambda, \Lambda}^{-}(M-N) - \delta \vert u-v\vert   \hspace{-.2cm} & \leq & \hspace{-.2cm} G(M,u,x) - G(N,v,x)  & \hspace{-.2cm}\\
      & \leq & \hspace{-.2cm}\mathcal{M}_{\lambda, \Lambda}^{+} (M-N) + \delta \vert u-v\vert,\end{array} 
    \end{equation*}\\

\begin{theorem}\label{T2}
   Let $\Omega \subset  \mathbb{R}^n$, $n \geq 2$, and assume $(H_0)$ - $(H_2)$ and $(D_F)$. Let $u$ a $L^p-$viscosity solution to 
   \begin{equation*}
    \left \{ \begin{array}{rclrl} 
    G(D^2u,u,x) \hspace{-.2cm} & \leq & \hspace{-.2cm} f(u) & \mbox{ in }  & \hspace{-.2cm} \Omega \\
     u(x) \hspace{-.2cm} & \geq & \hspace{-.2cm} 0 & \mbox{ in } & \hspace{-.2cm} \Omega\\
    u(x_{0}) \hspace{-.2cm} & = & \hspace{-.2cm} 0 & \mbox{ for some } & \hspace{-.2cm} x_0 \in \Omega. \end{array} \right.
\end{equation*}
with $f$ as in Theorem \ref{T3}. 

Then $u \equiv 0$.
\end{theorem}

First let's simplify the problem. The same argument from Theorem \ref{T3}, shows that we can assume $\Omega = \B_1$ and take $x_0 = 0$. By condition $(D_F)$, its sufficient to prove the result for the Pucci Operator:
    \begin{equation*}
    \left \{ \begin{array}{rclrl} 
    \mathcal{M}_{\lambda, \Lambda}^{-}(D^2u) - c(x)u(x) \hspace{-.2cm} & \leq & \hspace{-.2cm} f(u) & \mbox{ in }  & \hspace{-.2cm} \B_1 \\
     u(x) \hspace{-.2cm} & \geq & \hspace{-.2cm} 0 & \mbox{ in } & \hspace{-.2cm} \B_1\\
    u(0) \hspace{-.2cm} & = & \hspace{-.2cm} 0 &  \end{array} \right.
\end{equation*}

\begin{proposition}\label{Prop2}
   Assume the same hyphotesis of Theorem \ref{T2}. Let $u$ some solution of  
     \begin{equation*}
    \left \{ \begin{array}{rclrl} 
    \mathcal{M}_{\lambda, \Lambda}^{-}(D^2u) - c(x)u(x) \hspace{-.2cm} & \leq & \hspace{-.2cm} f(u) & \mbox{ in }  & \hspace{-.2cm} \B_1 \\
     u(x) \hspace{-.2cm} & \geq & \hspace{-.2cm} 0 & \mbox{ in } & \hspace{-.2cm} \B_1\\
    u(0) \hspace{-.2cm} & = & \hspace{-.2cm} 0 &  \end{array} \right.
\end{equation*}
Then we can find coeficients $\tilde{b}, \tilde{c} \in C(\B_{1/2})$, constants $K,C > 0$  and a solution $w \in C_{loc}^{1,\alpha}(\B_{1/2})$ of
\begin{equation*}
    \left \{ \begin{array}{rclrl} 
     \mathcal{M}_{\lambda, \Lambda}^{-}(D^2u) + \tilde{b}(y) \vert D w \vert - \tilde{c} w \hspace{-.2cm} & \leq & K f(C w)  & \mbox{ in }  & \hspace{-.2cm} \B_{1/2} \\
     w \hspace{-.2cm} & \geq & \hspace{-.2cm} 0 & \mbox{ in } & \hspace{-.2cm} \B_{1/2}\\
    w(0) \hspace{-.2cm} & = & \hspace{-.2cm} 0  \end{array} \right.
\end{equation*}
such that $w=0$ if and only if $u = 0$ in $\B_{1/2}$.
\end{proposition}
\begin{proof}
   First, observe that we can replace $c$ by your negative part $c^-$. Instead, how we know that $c = c^{+} - c^{-}$, then 
\begin{equation*}
    \mathcal{M}_{\lambda, \Lambda}^{-} (D^2 v(y)) - c^{-}(y)v(y) \leq \mathcal{M}_{\lambda, \Lambda}^{-} (D^2 v(y)) -c(y)v(y) \leq f(v(y)).
\end{equation*} 

Now, consider $\lambda_1^+ = \lambda_1^+(\mathcal{M}_{\lambda, \Lambda}^{-} (D^2 \cdot) - c^{-}\cdot, \B_1) > 0$ and $\varphi_1^+ \in W^{2,n}(\B_1)$ which are the eigenvalue and eigenfunction given by Theorem \ref{EAU2}. This means that $\varphi_1^+$ solves 
\begin{equation*}
     \left \{ \begin{array}{rclrl} 
   \mathcal{M}_{\lambda, \Lambda}^{-} (D^2 \varphi_1^+) - c^-(x) \varphi_1^+ \hspace{-.2cm} & = & \hspace{-.2cm} -\lambda_1^+ \varphi_1^+ & \mbox{ in }  & \hspace{-.2cm} \B_1 \\
     \varphi_1^+ \hspace{-.2cm} & > & \hspace{-.2cm} 0 & \mbox{ in } & \hspace{-.2cm} \B_1\\
    \varphi_1^+ \hspace{-.2cm} & = & \hspace{-.2cm} 0 & \mbox{ on } & \hspace{-.2cm} \partial \B_1 \end{array} \right.
\end{equation*}

The positivity of the first eigenvalue follows from the ABP estimates, wich imply the Comparision Principle, and we can apply the Theorem $1.5$ of \cite{AutV2}.

With this, we have $\inf_{\B_{1/2}}\varphi_1^+ > 0$, implying that the function $w : = \displaystyle\frac{v}{\varphi_1^+}$ is well-posed. Furthermore, with the help of Lemma \ref{LEM} and the positivity of $\varphi_1^+$ and $w$, the function $w$ satisfies:
\begin{equation*}
\begin{array}{rclrl}
 \mathcal{M}_{\lambda, \Lambda}^{-}(D^2(w \varphi_1^+)) - c^{-}(y)w \varphi_1^+  & \hspace{-.2cm} = & \hspace{-.2cm} \mathcal{M}_{\lambda, \Lambda}^{-}\left((D^2 w) \varphi_1^+ +  Dw\otimes D \varphi_1^+ + D \varphi_1^+ \otimes Dw +(D^2\varphi_1^+) w\right) \vspace{.1cm}\\ 
 &  & \hspace{-.1cm}- c^{-}(y)w \varphi_1^+ \vspace{.1cm} \\ 
 & \hspace{-.2cm} \geq & \hspace{-.2cm} \mathcal{M}_{\lambda, \Lambda}^{-}\left(D^2 w\right) \varphi_1^+ +  \mathcal{M}_{\lambda, \Lambda}^{-}\left( Dw\otimes D \varphi_1^+ \right) +  \mathcal{M}_{\lambda, \Lambda}^{-}\left(D \varphi_1^+ \otimes Dw \right) \vspace{.1cm} \\
 & & + \mathcal{M}_{\lambda, \Lambda}^{-}\left(D^2\varphi_1^+ \right) w - c^{-}(y)\varphi_1^+ w \vspace{.1cm}\\
 & \hspace{-.2cm} \geq & \hspace{-.2cm}\mathcal{M}_{\lambda, \Lambda}^{-}\left(D^2 w\right) \varphi_1^+ -2C(n,\lambda,\Lambda) \vert D\varphi_1^+\vert \vert D w \vert - w \lambda_1^+ \varphi_1^+. 
\end{array}
\end{equation*}
Therefore, $w$ is a solution for 
\begin{equation*}
    \left \{ \begin{array}{rclrl} 
     \mathcal{M}_{\lambda, \Lambda}^{-}\left(D^2 w\right) + \tilde{b}(y) \vert D w \vert  - \tilde{c}(y) w \hspace{-.2cm} & \leq & \hspace{-.2cm} K f(\varphi_1 w) \leq K f(C w)  & \mbox{ in }  & \hspace{-.2cm} \B_{1/2} \\
     w \hspace{-.2cm} & \geq & \hspace{-.2cm} 0 & \mbox{ in } & \hspace{-.2cm} \B_{1/2}\\
    w(0) \hspace{-.2cm} & = & \hspace{-.2cm} 0 &  \end{array} \right.
\end{equation*}
with:
\begin{equation*}
    K = \displaystyle(\inf_{\tilde{G}}\varphi_1)^{-1} > 0, \hspace{1mm} C = \|\varphi_1\|_{W^{2,n}(\tilde{G})}, \hspace{1mm} \tilde{b}(y) = -\frac{1}{\varphi(y)}2C(n,\lambda,\Lambda) \vert D\varphi_1^+(y)\vert \hspace{.1cm} \textrm{ and } \hspace{.1cm} \tilde{c}(y) = \lambda_1^+.
\end{equation*}
\end{proof}
Now, we can prove the last result of this work and it is imediate:\\

\textit{Proof of Theorem \ref{T2}.} This theorem is proved using a splitting argument for connected sets: we know that the set $\{u=0\}$ is closed, and we need to show that it is open. Since $\B_1$ is connected, we proceed by showing that there exists an open set $\tilde{B} \subset \B_1$ such that $\tilde{G} \subset \{u=0\}$.

By Proposition \ref{Prop2}, we obtain a solution $w$ defined in $\B_{1/2}$ such that $u(y) = \varphi_1(y) w(y)$ for $y \in \B_{1/2}$. We now need to show that $w \equiv 0$ on $\B_{1/2}$. This follows immediately from applying the  (VSMP) in the bounded setting proved in Theorem \ref{VBounded}. Note that we can, in fact, apply this result, because all the coefficients $ \tilde{b}$ and $\tilde{c}$ are bounded (contínuous), and $y_0$ is such that $w(y) = 0$. The remaining task is to show that $\tilde{f}(y) = f(Cy)$ satisfies the hypotheses on $f$. Let
\begin{equation*}
    \Tilde{F}(s) = \int_{0}^{s} f(\|\varphi_1\| t)dt
\end{equation*}
We need to show that
\begin{equation}\label{33}
    \int_{0}(\Tilde{F}(s))^{-1/2}ds = \infty.
\end{equation}
By the method of integration we have that
\begin{equation*}
    \Tilde{F}(s) = \int_{0}^{s} f(\|\varphi_1\| t)dt = \|\varphi_1\|^{-1}\int_{0}^{s\|\varphi_1\|}f(y)dy.
\end{equation*}
The hypothesis in equation $(\ref{33})$ tells us that the primitive tends to infinity as it approaches zero, and the change of variables affects only the upper part of the integral, multiplying it by a constant. This together with the original hypotesis on $f$ implies exactly what we intended.
\begin{center}
    \hspace{14.6cm}$\square$
\end{center}

 Now, we are ready to prove the main Theorem.\\
 
\textit{Proof of Theorem \ref{T3}}. After the reductions and applying the Proposition \ref{Prop1}, we just need to apply the Theorem \ref{T2} with $G = \mathcal{L}_{ND}$ and get the result.

\begin{center}
    \hspace{14.6cm}$\square$
\end{center}


We finalize this work with a usefull and interesting remark: the result still remains open for fully nonlinear operators that have dependence on the gradient. The Zvonkin transform, when applied to the Pucci operators, introduces a new lower-order term that has no definite sign, making it more difficult to work with the infimum.

\newpage
\bigskip


\bibliographystyle{amsplain, amsalpha}


\end{document}